\documentclass[12pt]{amsart}
\usepackage{amssymb}
\usepackage{graphicx}

\newcommand{\N}{{\mathbb N}}
\newcommand{\C}{{\mathbb C}}

\newcommand{\wand}{wandering domain}

\newcommand{\tef}{transcendental entire function}

\newcommand{\sconn}{simply connected}
\newcommand{\mconn}{multiply connected}

\theoremstyle{plain}
\newtheorem{theorem}{Theorem}[section]
\newtheorem{corollary}{Corollary}[section]

\newtheorem{lemma}{Lemma}[section]
\theoremstyle{definition}

\theoremstyle{remark}

\theoremstyle{problem}

\theoremstyle{example}

\begin{document}

\title[Escaping Fatou components]{Boundaries of escaping Fatou components}

\author{P. J. Rippon}
\address{Department of Mathematics and Statistics\\
The Open University \\
Walton Hall\\
Milton Keynes MK7 6AA\\
UK}
\email{p.j.rippon@open.ac.uk}

\author{G. M. Stallard}
\address{Department of Mathematics and Statistics\\
The Open University \\
Walton Hall\\
Milton Keynes MK7 6AA\\
UK}
\email{g.m.stallard@open.ac.uk}



\subjclass{30D05, 37F10.\newline\hspace*{.32cm} Both authors are
supported by EPSRC grant EP/H006591/1.}


\begin{abstract}
Let $f$ be a {\tef} and $U$ be a Fatou component of $f$. We show
that if $U$ is an escaping {\wand} of $f$, then most boundary
points of $U$ (in the sense of harmonic measure) are also
escaping. In the other direction we show that if enough boundary
points of $U$ are escaping, then $U$ is an escaping Fatou
component. Some applications of these results are given; for
example, if $I(f)$ is the escaping set of
$f$, then $I(f)\cup\{\infty\}$ is connected.
\end{abstract}

\maketitle

\section{Introduction}
\setcounter{equation}{0}
Let $f$ be a {\tef} and denote by $f^{n},\,n=0,1,2,\ldots\,$, the
$n$th iterate of~$f$. The {\it Fatou set} $F(f)$ is defined to be
the set of points $z \in \C$ such that $(f^{n})_{n \in \N}$ forms
a normal family in some neighborhood of $z$.  The components of
$F(f)$ are called {\it Fatou components}. The complement of
$F(f)$ is called the {\it Julia set} $J(f)$. An introduction to
the properties of these sets can be found in~\cite{wB93}.

The set $F(f)$ is completely invariant, so for any component $U$
of $F(f)$ there exists, for each $n=0,1,2,\ldots\,$, a component
of $F(f)$, which we call $U_n$, such that $f^n(U) \subset U_n$.
If, for some $p\ge 1$, we have $U_p =U_0= U$, then we say that
$U$ is a periodic component of {\it period} $p$, assuming $p$ to
be minimal. There are then four possible types of periodic
components; see \cite[Theorem~6]{wB93}. If $U_n\ne U_m$ for $m\ne
n$, then we say that $U$ is a {\it wandering domain} of $f$.

The {\it escaping set}
\[
 I(f) = \{z: f^n(z) \to \infty \mbox{ as } n \to \infty \}
\]
was first studied for a general transcendental entire function
$f$ by Eremenko~\cite{Ere89}. He proved that
\begin{equation}\label{erem}
I(f)\cap J(f) \neq \emptyset\quad\text{and}\quad \partial
I(f)=J(f),
\end{equation}
and also that $\overline{I(f)}$ has no bounded
components. Eremenko remarked that it is plausible that all
the components of $I(f)$ are unbounded, a statement now known as
{\it Eremenko's conjecture} that remains open in spite of much
work on it and many partial results.

Any Fatou component that meets $I(f)$ must lie in $I(f)$ by
normality; we call such components {\it escaping Fatou
components}. Escaping wandering domains can be bounded or unbounded, and escaping periodic Fatou components are all
unbounded; the latter are called {\it Baker domains}. This paper
gives results about the relationship between an escaping Fatou
component and its boundary, and describes some consequences of
these results.

For the function $f(z)=z+1+e^{-z}$, studied by Fatou in
\cite{pF26}, the set $F(f)$ is a completely invariant Baker
domain, whose boundary is $J(f)$. So in this example an escaping
Fatou component has many boundary points that are {\it not} in
$I(f)$. It is natural to ask whether every escaping Fatou
component of a {\tef} must have at least one boundary point in
$I(f)$. We have the following partial result.

\begin{theorem}\label{Thm1.1}
Let $f$ be a {\tef} and let $U$ be a {\wand} of $f$ such that
$U\subset I(f)$. Then
\[\partial U\cap I(f) \ne \emptyset.\]
Moreover, the set $\partial U\cap I(f)^c$ has zero harmonic measure
relative to $U$.
\end{theorem}

The proof of Theorem~\ref{Thm1.1} can be adapted to show that
$\partial U\cap I(f) \ne \emptyset$ for many Baker domains (see
Section~2, Remark~2), but it remains open whether this conclusion
holds whenever $U$ is a Baker domain.

In the other direction to Theorem~\ref{Thm1.1}, we can ask
whether a Fatou component must be escaping if a large enough
subset of its boundary is escaping. We have the following result in this direction.
\begin{theorem}\label{Thm1.2}
Let $f$ be a {\tef} and let $U$ be a Fatou component of $f$.
\begin{itemize}
\item[(a)] If \,$\partial U\cap I(f)$ has positive harmonic
measure relative to $U$, then $U\subset I(f)$.
 \item[(b)]
If \,$\partial U\cap A(f)$ has positive harmonic measure relative
to $U$, then $U\subset A(f)$.
\end{itemize}
\end{theorem}
In particular, if $U$ is a Fatou component of $f$ such that
$\partial U\subset I(f)$, then $U\subset I(f)$, and a similar
result holds for $A(f)$.

Here $A(f)$ is the {\it fast escaping set},
introduced by Bergweiler and Hinkkanen in~\cite{BH99}, which can be
defined as follows (see~\cite{RS09}):
\[
A(f) = \{z: \mbox{there exists } L \in \N \mbox{ such that }
|f^{n+L}(z)| \ge M^n(R), \mbox{ for } n\in \N\}.
\]
Here \(M(r)= M(r,f) = \max_{|z|=r} |f(z)| \) and $R>0$ is such
that $M(r)>r$ for $r\ge R$ or equivalently such that $M^n(R)\to
\infty$ as $n\to\infty$.

Many of the properties of $A(f)$ are stronger than those of $I(f)$ (see Section~4).
For example, $A(f)$ has the following stronger property
than that given for $I(f)$ in Theorem~\ref{Thm1.1} (see~\cite[Theorem~1.2]{RS09}):
\begin{equation}\label{eqn1.1}
\text{if }U \text{ is a Fatou component of }f\text{ that meets
}A(f), \text{ then }\overline{U}\subset A(f).
\end{equation}

Note that it is possible for a Fatou component to lie in $A(f)$;
for example,
\begin{equation}\label{eqn1.2}
\text{if }U \text{ is a multiply connected Fatou component of }f,\text{ then }\overline{U}\subset A(f),
\end{equation}
(see~\cite{RS05}) and Bergweiler has constructed a
{\tef} with both {\sconn} and {\mconn} Fatou components
in $A(f)$ (see~\cite{wB09}). All Fatou components in $A(f)$ are wandering domains (see \cite{BH99}).

We prove Theorems~\ref{Thm1.1} and~\ref{Thm1.2} in Sections~2
and~3, respectively. In Section~4 we give two consequences of
these theorems, one of which states that if $f$ is a {\tef}, then
$I(f)\cup\{\infty\}$ is connected, and the other of which is an
improved version of a result from~\cite{RS07} which gives a
sufficient condition for $I(f)$ to be connected; the proof we
give of the latter result was inspired by an unpublished idea of
Professor Noel Baker. This section also includes a short proof
that $I(f)\ne\emptyset$, which may be of independent interest.

Section~5 is concerned with the components of $I(f)$; for
example, we give a new sufficient condition for the components of
$I(f)$ to be unbounded and prove that various sets, such as
$I(f)$ and $A(f)$, and their complements, are either connected or
have infinitely many components.

{\it Acknowledgements}\quad We thank Walter Bergweiler, Dan Nicks
and Dave Sixsmith for their useful comments, and Lasse Rempe for
a nice observation about the statement of Theorem~\ref{Thm4.1}.

\section{Proof of Theorem~\ref{Thm1.1}}
\setcounter{equation}{0} First we recall that for a domain $G$
and a set $E\subset\partial G$, the {\it harmonic measure} of $E$
at $z$ relative to $G$, denoted by $\omega(z,E,G)$, is the
solution of the Dirichlet problem in $G$ (found by using the
Perron method) with boundary values given by the characteristic
function $\chi_E$. See \cite{GM05} or \cite{tR95}, for example,
for the solution of the Dirichlet problem and the properties of
harmonic measure.

In any {\sconn} domain $G$ in $\C$ all points $\zeta\in
\partial G$ are {\it regular} for the Dirichlet problem (see~\cite{tR95}); that
is, if $\phi$ is any real-valued function on $\partial G$ which
is continuous at $\zeta\in\partial G$ and $H_{\phi}$ is the
solution of the Dirichlet problem in $G$ with boundary
values~$\phi$, then $H_{\phi}(z)\to \phi(\zeta)$ as $z\to \zeta$.
\begin{proof}[Proof of Theorem~1.1]
Suppose that $U$ is an escaping wandering domain of $f$. Fix
$z_0\in U$ and, for $n\ge 0$, let $z_n=f^n(z_0)$ and $U_n$ be the
Fatou component that contains $f^n(U)$. In view of (\ref{eqn1.2}), we can assume that each $U_n$ is {\sconn}.

Let $R>0$ and put
\[
B_n=B_n(R)=\{z\in \partial U: |f^n(z)|\le R\},\quad n\in\N.
\]
To prove that $\partial U\cap I(f)^c$ has zero harmonic measure relative to $U$, it is sufficient to show that, for each $R>0$,
\begin{equation}\label{eq2.1}
\bigcap_{m\ge 1}\bigcup_{n\ge m}B_n=\{z\in\partial U:|f^n(z)|\le R \;\text{for infinitely many}\;n\}
\end{equation}
has harmonic measure~0 relative to $U$.

Define $\Delta=\{z:|z|>R\}\cup\{\infty\}$, considered as a disc
in the Riemann sphere $\hat{\C}$, and choose $N=N(R)$ such that
$|z_n|>2R$ for $n\ge N$. Then, for $n\ge N$, define
\[
E_n=\partial U_n\cap \{z:|z|\le R\},
\]
\[
V_n\;\text{to be the component of}\;U_n\cap \Delta\;\text{that contains}\;z_n,
\]
\[
F_n=\partial V_n\cap \{z:|z|= R\};
\]
see Figure~1.
\begin{figure}[htb]
\begin{center}
\includegraphics[width=9.0cm]{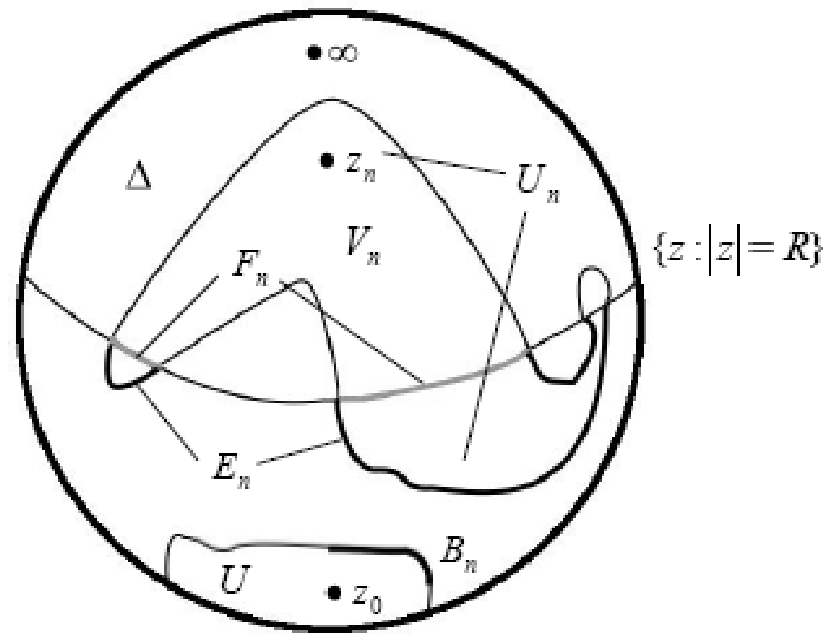}
\caption{The sets $U$, $B_n$, $\Delta$, $U_n$, $V_n$, $E_n$ and
$F_n$}
\end{center}
\end{figure}
We have
\begin{equation}\label{eq2.2a}
\omega(z,E_n, U_n)\le \omega(z,F_n,\Delta), \quad\text{for}\;
z\in V_n, n\ge N.
\end{equation}
This inequality is a special case of Carleman's domain extension principle, and it is proved as follows. The function
\[
u(z)=\left\{
 \begin{array}{ll}
 0,&z\in \Delta\setminus V_n,\\
 \omega(z,E_n,U_n),&z\in V_n,
 \end{array}
 \right.
\]
is subharmonic in $\Delta$ because $\omega(z,E_n,U_n)$ has boundary value $0$ at all points of $\partial V_n\cap\Delta$ (since $V_n$ is {\sconn} and hence regular for the Dirichlet problem). Also,
\[
\limsup_{z\to\zeta} u(z)\le \left\{
 \begin{array}{ll}
 0,&\text{for }\zeta\in \partial \Delta\setminus F_n,\\
 1,&\text{for }\zeta\in F_n,
 \end{array}
 \right.
\]
because $F_n$ is a closed subset of $\partial \Delta$. So (\ref{eq2.2a}) holds by the maximum principle for subharmonic functions (see \cite[Theorem~2.3.1]{tR95}).

By conformally mapping $\Delta$ onto $\{z:|z|<1/R\}$ and applying Harnack's inequality (see \cite[Theorem~1.3.1]{tR95}) in that disc, we obtain
\begin{equation}\label{eq2.2b}
\omega(z_n,F_n,\Delta)\le 3\,\omega(\infty,F_n,\Delta), \quad\text{for}\; n\ge N,
\end{equation}
as $|z_n|>2R$ for $n\ge N$. Since $U$ is a {\wand}, the sets $F_n$ are disjoint, so we deduce from \eqref{eq2.2a} and \eqref{eq2.2b} that
\begin{equation}\label{eq2.2}
\sum_{n\ge N}\omega(z_n,E_n, U_n)\le 3\sum_{n\ge
N}\omega(\infty,F_n,\Delta)\le 3\,\omega(\infty,\partial \Delta,
\Delta)= 3.
\end{equation}
Now $f^n(U)\subset U_n$ and $f^n(B_n)\subset E_n$, so
\[
\omega(z,B_n, U)\le \omega(f^n(z),E_n,U_n), \quad\text{for}\; z\in U, n\ge N,
\]
by \cite[Theorem~4.3.8]{tR95}. In particular,
\[
\omega(z_0,B_n, U)\le \omega(z_n,E_n,U_n), \quad\text{for}\; n\ge N.
\]
Hence, by~(\ref{eq2.2}),
\begin{equation}\label{eq2.3}
\sum_{n\ge N}\omega(z_0,B_n, U)\le \sum_{n\ge
N}\omega(z_n,E_n,U_n)<\infty,
\end{equation}
so
\[
\omega(z_0,\textstyle\bigcup_{n\ge m}B_n,U)\le \sum_{n\ge m}\omega(z_0,B_n,U)\to 0\quad \text{as } m\to\infty,
\]
which gives (\ref{eq2.1}).
\end{proof}
{\it Remarks}

1. It is natural to ask if there is an example of an entire
function with an escaping wandering domain that has at least one
non-escaping boundary point.

2. The proof of Theorem~\ref{Thm1.1} can be adapted to show that,
for many Baker domains $U$, the set $\partial U\cap I(f)^c$ has
zero harmonic measure relative to $U$. For example, suppose that
$U$ is an invariant Baker domain in which there is an orbit
$z_n=f^n(z_0)$, $n\in\N$, such that
\begin{equation}\label{eq2.4}
|z_{n+1}|\ge k|z_n|,\quad\text{for }n\in\N,
\end{equation}
where $k>1$. This is the case for several types of Baker domains;
see~\cite[Section~2]{pR08} for examples. Then, with the notation
from the proof of Theorem~\ref{Thm1.1}, we have $U_n=U$ and
$E_n=\partial U\cap\{z:|z|\le R\}=E$, say, for $n\in\N$, and
\begin{equation}\label{eq2.5}
\omega(z,E,U)\le \frac{C}{|z|^{1/2}},\quad\text{for }z\in U\cap \Delta,
\end{equation}
where $C>0$ is an absolute constant. Since $U$ is {\sconn}, this
last harmonic measure estimate can be obtained from the Beurling
projection theorem (see \cite[Theorem~9.2, page~105]{GM05}), by
exchanging the roles of~$0$ and~$\infty$. By (\ref{eq2.4}) and
(\ref{eq2.5}), we deduce that (\ref{eq2.3}) holds, so the
conclusion of Theorem~\ref{Thm1.1} also holds in this situation.

\section{Proof of Theorem~\ref{Thm1.2}}
\setcounter{equation}{0} %
The proof of Theorem~\ref{Thm1.2} uses the following version of the maximum principle
for subharmonic functions (see \cite[page~102]{GM05} or
\cite[Theorem~3.6.9]{tR95}).
\begin{lemma}\label{lem3.1}
Let $G$ be a domain in $\C$ and let $u$ be a subharmonic function
in~$G$ which is bounded above. If $\partial G$ is not a polar set
and
\[
\limsup_{z\to\zeta}u(z)\le 0,\quad\text{for }\zeta\in\partial
G\setminus E,
\]
where $E$ is a polar subset of $\partial G$, then $u\le 0$ in $G$.
\end{lemma}
Polar sets are defined in \cite{GM05} and \cite{tR95} -- here we just need the facts that if a set contains
a continuum, then it is not polar whereas a finite set is polar.
\begin{proof}[Proof of Theorem~\ref{Thm1.2}]
In view of (\ref{eqn1.2}), we can assume throughout the proof
that~$U$ is {\sconn}; in particular, $\partial U$ is not a polar
set and all points of $\partial U$ are regular for the Dirichlet
problem.

{\it Part~(a).}\quad We consider two cases. First we suppose that
$F(f)$ is disconnected. In this case we claim that there is a
closed disc which lies outside $\bigcup_{n\ge 0}f^n(U)$.
This is clearly true if $U$ is forwards invariant under $f$. If
$U$ is not forwards invariant, then there is a Fatou component
$V$ of $f$ different from $U$ such that $f(V)\subset U$.
Moreover, we claim that we can choose $V$ to be different from all $f^n(U)$,
$n\ge 0$. For otherwise, $U$ is periodic under $f$, with period
$p>1$, and the~$p$\, Fatou components of $f$ in the orbit of $U$
are each completely invariant under~$f^p$. This contradicts a
theorem of Baker~\cite{inB70} which states that a {\tef} can have at most one completely invariant Fatou component.

By making an affine change of variables, we can assume that
$\{z:|z|\le 1\}$ lies outside $\bigcup_{n\ge 0}f^n(U)$. Thus the functions
\[
u_n(z)=\log|f^n(z)|,\quad  z\in \overline{U},\; n\in\N,
\]
are positive harmonic in $U$ and continuous in
$\overline{U}$. By hypothesis there exists a set $E\subset
\partial U$ such that the harmonic measure $\omega(z,E,U)>0$, for
$z\in U$, and
\[
u_n(\zeta)\to\infty\;\text{ as }n\to \infty,\quad \text{for }\zeta\in E.
\]

Since $\omega(z,\,.\,,U)$ is a positive finite Borel measure on
$\partial U$, we can assume by Egorov's theorem (see \cite{wR70},
for example) that
\[
u_n\to\infty\;\text{ as }n\to \infty,\quad \text{uniformly on }
E,
\]
and we can assume in addition that $E$ is closed. Thus if $C>0$ is
given, then there exists $N=N(C)$ such that
\begin{equation}\label{lower}
u_n(\zeta)\ge C,\quad\text{for } \zeta\in E,\; n\ge N.
\end{equation}
Now consider the bounded harmonic functions
 \[v_n(z)=C\omega(z,E,U)-u_n(z),\quad z\in U,\; n\ge N.\]
 Then, for $n\ge N$,
 \[
 \limsup_{z\to \zeta} v_n(z)\le
 \left\{
 \begin{array}{ll}
 C-u_n(\zeta)\le 0,&\text{for }\zeta\in E,\\
 \displaystyle\lim_{z\to \zeta} C\omega(z,E,U)-u_n(\zeta)\le 0,&\text{for }\zeta\in \partial U\setminus E.
 \end{array}
 \right.
 \]
The first inequality follows from~(\ref{lower}) and the fact that $\omega(z,E,U)\le 1$ for $z\in
U$. The second inequality holds because
$u_n(z)>0$ for $z\in U$ and $\lim_{z\to \zeta}\omega(z,E,U)=0$,
for $\zeta\in\partial U\setminus E$, since all points of
$\partial U$ are regular and $E$ is closed, and $\omega(z,E,U)$
is the solution of the Dirichlet problem in $U$ with boundary values
$\chi_E$.

Also, in the case that $U$ is unbounded, the point at infinity
acts as a polar set in the boundary of $U$ (see \cite[p.~85 and
Corollary~3.2.5]{tR95}).

It follows by Lemma~\ref{lem3.1} that $v_n(z)\le 0$, for $z\in U$.
Hence
\[
u_n(z)=\log|f^n(z)|\ge C\omega(z,E,U),\quad\text{for }z\in U,\;
n\ge N.
\]
Since $C>0$ is arbitrary, it follows that
\[u_n(z)=\log|f^n(z)|\to\infty\;\text{ as } n\to\infty,\quad\text{for }z\in U,\]
so $U\subset I(f)$.

In the second case, we suppose that $U=F(f)$ is connected, from
which it follows that $f^n(U)\subset U$, for $n\in\N$. In this
case, $U$ is {\sconn} and unbounded so we can define
\begin{equation}\label{psi}
\psi(z)=k\sqrt{z-a}+b,\quad z\in U,
\end{equation}
where $a\in J(f)$, $\psi$ is conformal on $U$, and $k>0$ and
$b\in\C$ are chosen so that $\psi(U)\cap \{z:|z|\le
1\}=\emptyset$.
 Note that for $z\in U$ we have $\psi(z) \to \infty$ if and only if $z\to\infty$. Then put
\[
v_n(z)=\log|\psi(f^n(z))|,\quad z\in\overline{U},\; n\in\N.
\]
Each $v_n$ is positive harmonic in $U$ and continuous in
$\overline{U}$, and we find by using Egorov's theorem again that
there exists a closed set $E\subset \partial U$ such that
$\omega(z,E,U)>0$ for $z\in U$ and
\[
v_n\to\infty\;\text{ as }n\to \infty,\quad \text{uniformly on }
E.
\]
As in the first case, we can use Lemma~\ref{lem3.1} to deduce
that
\[
v_n(z)=\log|\psi(f^n(z))|\to \infty \;\text{ as }
n\to\infty,\quad\text{for }z\in U.
\]
Hence $f^n(z)\to \infty$ as $n\to\infty$ for
$z\in U$, as required.

{\it Part~(b)}\quad This proof is similar to that of part~(a).
Suppose that $\partial U\cap A(f)$ has positive harmonic measure
relative to $U$. By the definition of $A(f)$, we have
\[A(f)=\bigcup_{L\in\N}A_R^{-L}(f),\]
where $R>1$ is so large that $R_n=M^n(R)\to \infty$ as
$n\to\infty$ and
\[A_R^{-L}(f)=\{z:|f^{n+L}(z)|\ge R_n,\;\text{for }n\in\N\},\quad L\in\N.\]
Thus $\partial U\cap A_R^{-L}(f)$ has positive harmonic measure
relative to $U$ for some $L\in\N$. Hence there exists a closed
set $E\subset\partial U$ of positive harmonic measure relative
to~$U$ such that
\begin{equation}\label{eq3.2}
|f^{n+L}(\zeta)|\ge R_n,\quad \text{for } \zeta\in E,\; n\in\N.
\end{equation}
As in the proof of part~(a), in the `$F(f)$ disconnected' case we
deduce by Lemma~\ref{lem3.1} that
\[
\log|f^{n+L}(z)|\ge \left(\log R_n\right)\,\omega(z,E,U),\quad
\text{for } z\in U,\; n\in\N,
\]
so
\[
|f^{n+L}(z)|\ge R_n^{\omega(z,E,U)},\quad \text{for } z\in U,\;
n\in\N.
\]
For any fixed $z\in U$, we can choose $N$ so large that
\[
R_n^{\omega(z,E,U)}\ge R_{n-1},\quad\text{for } n\ge N,
\]
because $\log M(r)/\log r\to \infty$ as $r\to\infty$. So for this
$z$ we have
\[
|f^{n+L}(z)|\ge R_{n-1},\quad \text{for } n\ge N,
\]
and, moreover,
\[
|f^{n+L+1}(z)|\ge R_{n},\quad \text{for } n\in\N,
\]
because  $|f(z)|\ge M(r)$ implies that $|z|\ge r$ for $r>0$. Thus
$z\in A_R^{-L-1}(f)$.

The proof in the `$F(f)$ connected' case is similar, using the
fact that, for $|z|$ large enough, the function $\psi$ defined in (\ref{psi}) satisfies
\begin{equation}\label{psiestimate}
\frac{k}{2}|z|^{1/2}\le|\psi(z)|\le 2k|z|^{1/2}.
\end{equation}
It follows from (\ref{eq3.2}) and (\ref{psiestimate}) that, for some $N\in\N$, we have
\[
|\psi(f^{n+L}(\zeta))|\ge (k/2)R_n^{1/2},\quad \text{for }
\zeta\in E,\; n\ge N,
\]
so, using Lemma~\ref{lem3.1} as before, we deduce that
\[
\log|\psi(f^{n+L}(z))|\ge \left(\log
\left((k/2)R_n^{1/2}\right)\right)\,\omega(z,E,U),\quad \text{for
} z\in U,\; n\ge N.
\]
Hence
\[
|\psi(f^{n+L}(z))|\ge
\left((k/2)R_n^{1/2}\right)^{\omega(z,E,U)},\quad \text{for }
z\in U,\; n\ge N,
\]
so, by~(\ref{psiestimate}) again, there exists $c>0$ such that
\[
|f^{n+L}(z)|\ge cR_n^{\omega(z,E,U)},\quad \text{for } z\in U,\;
n\ge N.
\]
The rest of the proof is similar to that of the `$F(f)$
connected' case in part~(a).
\end{proof}
{\it Remarks}

1. If the Fatou component $U$ is bounded, then the proofs of both
parts of Theorem~\ref{Thm1.2} can be simplified considerably
since only the `$F(f)$ disconnected' case can occur, and
Lemma~\ref{lem3.1} can be replaced by the ordinary maximum
principle.

2. In Theorem~\ref{Thm1.2} part~(a), the condition that $\partial
U\cap I(f)$ has Hausdorff dimension~2 would not imply that
$U\subset I(f)$. For example, functions of the form $f(z)=\lambda
e^z$, $0<\lambda <1/e$, have a completely invariant attracting
basin~$U$ and the Hausdorff dimension of $\partial U\cap
I(f)=J(f) \cap I(f)$ is~2; see \cite{McM87}.


\section{Some consequences of Theorems~\ref{Thm1.1} and~\ref{Thm1.2}}
\setcounter{equation}{0}  In this section we prove two results
about $I(f)$ which follow from Theorems~\ref{Thm1.1}
and~\ref{Thm1.2}, respectively. We thank Lasse Rempe for pointing
out that part~(b) of Theorem~\ref{Thm4.1} can be stated in the
more interesting way given in part~(c).
\begin{theorem}\label{Thm4.1}
Let $f$ be a {\tef}. The following statements hold.
\begin{itemize}
\item[(a)] Any bounded component of $I(f)$ meets $J(f)$.
\item[(b)] If $G$ is a bounded simply connected domain and $G\cap
I(f)\ne \emptyset$, then $\partial G\cap I(f)\ne \emptyset$.
\item[(c)] $I(f)\cup \{\infty\}$ is connected.
\end{itemize}
\end{theorem}
To prove Theorem~\ref{Thm4.1} we need a lemma whose proof is based on the {\it blowing up property} of $J(f)$:
\begin{quote}
if $f$ is an entire function, $K$ is compact, $K\subset
\C\setminus E(f)$ and $V$ is an open neighbourhood of $z \in
J(f)$, then there exists $N\in\N$ such that $f^n(V)\supset
K$, for all $n\ge N$.
\end{quote}

Here $E(f)$ is the {\it exceptional set} of $f$, that is, the set of
points with a finite backwards orbit under $f$, which has at most
one point (see \cite{wB93}).
\begin{lemma}\label{lem5.1}
Let $f$ be a {\tef}. If $G$ is a bounded simply connected domain such that $G\cap J(f)\ne
\emptyset$, then $\partial G\cap I(f) \ne \emptyset$.
\end{lemma}
\begin{proof}
Since $G$ is a bounded domain, we can define $\alpha_n$, $n\in\N$,
to be the outer boundary component of $f^n(G)$. Then, by the
blowing up property of $J(f)$, we have dist\,$(\alpha_n,0)\to\infty$ as
$n\to\infty$. The compact sets
\[K_n=\{z\in\partial G: f^n(z)\in \alpha_n\}, \quad n\in\N,\]
form a nested sequence, since $\alpha_{n+1}\subset f(\alpha_n)$,
for $n\in\N$, so $K=\bigcap_{n\ge 1}K_n\ne\emptyset$. All points
of $K$ must lie in $I(f)$, so $\partial G\cap I(f)\ne \emptyset$, as required.
\end{proof}

\begin{proof}[Proof of Theorem~\ref{Thm4.1}]
Part~(a) follows immediately from Theorem~\ref{Thm1.1}. For if $I_0$ is a bounded component of $I(f)$ that does not meet $J(f)$, then $I_0$ must be a bounded Fatou component, and this contradicts the fact that such a component must have a boundary point in $I(f)$, by Theorem~\ref{Thm1.1}.

To prove part~(b) we suppose that $G$ is a bounded {\sconn}
domain that meets $I(f)$. By Lemma~\ref{lem5.1}, we can assume
that $G\subset F(f)$, so $G\subset U$ where $U$ is an escaping
Fatou component of~$f$. There are now two cases:
\begin{itemize}
 \item either $\partial G\cap U\ne \emptyset$, in which case
$\partial G\cap I(f)\ne \emptyset$;
 \item or $G=U$, in which case $U$ is an escaping wandering domain,
 so $\partial G\cap I(f)\ne \emptyset$, by Theorem~\ref{Thm1.1}.
\end{itemize}

To prove part~(c) we note that if $I(f)\cup \{\infty\}$ is
disconnected, then there exist disjoint open sets $H_1$ and $H_2$
in the Riemann sphere $\hat{\C}$ such that
\[I(f)\cup \{\infty\}\subset H_1\cup H_2\quad\text{and}\quad H_i\cap (I(f)\cup \{\infty\})\ne \emptyset,\;\text{for }i=1,2.\]
Without loss of generality $\infty\in H_2$ and $H_1$ is bounded.
 Since $I(f)$ meets $H_1$ it also meets $\partial H_1$ by
part~(b), and this gives a contradiction.

(Part~(c) clearly implies part~(b), so these two statements are equivalent.)
\end{proof}
{\it Remarks}

1. The proof above of Lemma~\ref{lem5.1} gives a surprisingly
simple argument to show that $I(f)$ is non-empty.

2. An alternative proof of Lemma~\ref{lem5.1} can be given, based
on the following properties of the fast escaping set $A(f)$
(see~\cite{BH99} and~\cite{RS05}):
\begin{equation}\label{eremA1}
\partial A(f)=J(f);
\end{equation}
\begin{equation}\label{eremA2}
\text{all the components of } A(f)\text{ are unbounded}.
\end{equation}

Since $G\cap J(f)\ne \emptyset$, we have $G\cap A(f)\ne\emptyset$, by (\ref{eremA1}). Hence $\partial G\cap A(f)\ne\emptyset$, by (\ref{eremA2}), which proves Lemma~\ref{lem5.1}. In fact, the set $K$ constructed in the proof of Lemma~\ref{lem5.1} is a subset of $A(f)$ by \cite[Corollary~2.5]{RS09}.

3. There is also a constructive proof of Theorem~\ref{Thm4.1}
part~(c). First prove Theorem~\ref{Thm4.1} part~(a) as above.
Then consider $A_0=A(f) \cup\{\infty\}$. By (\ref{eremA2}), $A_0$
is a union of connected sets in $\hat{\C}$, all containing
$\infty$, so $A_0$ is connected. Let $I_0$ be the component of
$I(f) \cup \{\infty\}$ that contains $A_0$. Any point of $I(f)$
either lies in an unbounded component of $I(f)$, and so is in
$I_0$, or it lies in a bounded component of $I(f)$, which must
contain a point of $J(f)$, by Theorem~\ref{Thm4.1} part~(a). By
(\ref{eremA1}), such a component of $I(f)$ must also be in $I_0$.
Hence $I(f) \cup\{\infty\}=I_0$.

Next we use Theorem~\ref{Thm1.2} to give a sufficient condition for $I(f)$ to be connected.
\begin{theorem}\label{Thm4.2}
Let $f$ be a {\tef} and suppose that there exists a bounded
domain $G$ such that
\begin{equation}\label{cond1}
\partial G\subset I(f)\quad\text{and}\quad G\cap I(f)^c\ne\emptyset.
\end{equation}
Then
\begin{itemize}
 \item[(a)] for each $n\in\N$, the outer boundary component
$\alpha_n$ of $f^n(G)$ is contained in $I(f)$, $\alpha_n\to
\infty$ as $n\to\infty$, and $\alpha_n$ surrounds~0 for
sufficiently large $n$;
 \item[(b)] $I(f)$ is connected.
\end{itemize}
\end{theorem}
The statement of Theorem~\ref{Thm4.2} is similar to that of
\cite[Theorem~2]{RS07}, but here we assume that $G$ satisfies (\ref{cond1}) whereas in \cite[Theorem~2]{RS07} we
assumed that
\begin{equation}\label{cond2}
\partial G\subset I(f)\quad\text{and}\quad G\cap
J(f)\ne\emptyset.
\end{equation}
In \cite{RS07}, we remarked that (\ref{cond1}) and
(\ref{cond2}) are equivalent, but we later realised that the
statement that (\ref{cond1}) implies (\ref{cond2}) is not
immediate. The following lemma shows that this
implication can be deduced from Theorem~\ref{Thm1.2}.
\begin{lemma}\label{Lem4.1}
Let $f$ be a {\tef}. Then the conditions (\ref{cond1}) and
(\ref{cond2}) are equivalent.
\end{lemma}
\begin{proof}
It is clear that if $G \cap J(f) \ne \emptyset$, then $G\cap
I(f)^c \ne \emptyset$, since $J(f)=\partial I(f)$, by (\ref{erem}). Thus
(\ref{cond2}) implies (\ref{cond1}).

To prove that (\ref{cond1}) implies (\ref{cond2}) we argue as
follows. If $\partial G\subset I(f)$ and $G \subset F(f)$, then
either $\partial G\cap F(f)\ne \emptyset$, in which case
$G\subset I(f)$ or $G$ is a component of the Fatou
set, so $G\subset I(f)$ by Theorem~\ref{Thm1.2} part~(a). Hence
(\ref{cond1}) implies (\ref{cond2}).
\end{proof}
Theorem~\ref{Thm4.2} follows immediately from Lemma~\ref{Lem4.1}
and \cite[Theorem~2]{RS07}. However, the proof below contains a
shorter argument to prove part~(b) than that given
in~\cite{RS07}. (In \cite{RS07} a constructive argument was used
whereas here we argue by contradiction.)
\begin{proof}[Proof of Theorem~4.2]
By Lemma~\ref{Lem4.1}, the hypothesis (\ref{cond1}) implies that
$G\cap J(f)\ne\emptyset$.

Since $\partial G\subset I(f)$, the outer boundary component
$\alpha_n$ of $f^n(G)$ is contained in $I(f)$, for each $n\in\N$.
By the blowing up property of $J(f)$, $\alpha_n\to \infty$ as
$n\to\infty$, and $\alpha_n$ surrounds~0 for sufficiently large
$n$. Thus part~(a) holds.

To prove part~(b) suppose that $I(f)$ is disconnected. Then there
exist disjoint open sets $H_1$ and $H_2$ in $\C$ such that
\[I(f)\subset H_1\cup H_2\quad\text{and}\quad H_i\cap I(f)\ne \emptyset,\; \text{for }i=1,2.\]
Without loss of generality one of $H_1$ and $H_2$ is bounded, say
$H_1$, for otherwise both $H_1$ and $H_2$ have unbounded boundary
components that do not meet $I(f)$, which contradicts part~(a).
We can also assume that $H_1$ is simply connected.

Fix $z_0\in H_1$ and take any $m\in\N$. Then there exists
$N_m\in\N$ such that
\begin{equation}\label{out} f^n(z_0)\;\text{lies outside}\;
\alpha_{m},\quad \text{for } n\ge N_m.
\end{equation}
 Now, for $n\in\N$,
\begin{equation}\label{boundary}
\partial f^n(H_1)\subset f^n(\partial H_1)\subset I(f)^c,
\end{equation}
because $\partial H_1\subset I(f)^c$. Also, $f^n(\partial
H_1)$ is connected because $\partial H_1$ is connected. Thus,
since $\alpha_m\subset I(f)$, it follows from (\ref{out}) and
(\ref{boundary}) that
\[\overline{f^n(H_1)}\;\text{lies outside}\; \alpha_{m},\quad \text{for } n\ge N_m.\]
Thus $\overline{f^n(H_1)}\to\infty$ as $n\to\infty$, so $\partial H_1\subset I(f)$, a contradiction. Hence $I(f)$ is connected.
\end{proof}

{\it Remarks}

1. The final part of this contradiction argument can be replaced
by an argument based on the use of Lemma~\ref{lem5.1} and
Theorem~\ref{Thm4.1}.

2. In Theorem~\ref{Thm4.2} the set $I(f)$ has the structure of a
`spider's web', a concept introduced in \cite{RS09}: we say that
a set~$E$ is an {\it (infinite) spider's web} if~$E$ is connected
and there exists a sequence $(G_n)$ of bounded simply connected
domains with $G_n\subset G_{n+1}$, $\partial G_n \subset E$, for
$n\in\N$, and
\[\bigcup _{n=1}^{\infty}G_n=\C.\]
Several examples of functions for which $I(f)$ is a spider's web,
and hence $I(f)$ is connected, were given
in~\cite[Section~6]{RS07}, though the name `spider's web' was not
used there. Many more examples of such functions are given in
\cite{RS09} and \cite{dS10}.

Condition (\ref{cond1}) can be described by saying that $I(f)$
{\it has a hole}, so Theorem~\ref{Thm4.2} can be stated in the following geometric way:
\begin{quote}
if $I(f)$ has a hole, then $I(f)$ is a spider's web.
\end{quote}

3.\; The proof of Theorem~\ref{Thm4.2} given above was inspired
by unpublished work of Professor Noel Baker. In early 2010 the
authors found a note written by Noel Baker some time after 1996 in an old diary that he
used for rough work. He stated there that if $f$ is a {\tef} with a {\mconn} Fatou component,
then $I(f)$ is connected. He gave no details of the proof, except
that it used the method of showing that $I(f)\ne\emptyset$ given
by Dom\'inguez in~\cite{pD98}, and there was a sketch which
suggested that his argument was by contradiction.

Noel Baker died in 2001 and as far as we know he never mentioned
this result to anyone. In \cite{RS05} we proved the result as a
corollary of the fact that all the components of $A(f)$ are
unbounded. Thinking about the argument he might have found led us
to the proof of Theorem~\ref{Thm4.2} above.

\section{Components of $I(f)$}
\setcounter{equation}{0} In this final section, we give various results
about the possible structures of the components of $I(f)$. The
first result, which follows from Theorem~\ref{Thm1.1}, gives a
new sufficient condition for all the components of $I(f)$ to be
unbounded.
\begin{theorem}\label{Thm3.1}
Let $f$ be a {\tef} and $E$ be a set such that
$J(f)\subset\overline{E}$. If $E$ is contained in the union of
finitely many components of $I(f)$, then
\begin{itemize}
 \item[(a)]\(I(f)\cap J(f)\) is contained in one component, $I_1$ say, of
$I(f)$;
 \item[(b)]all the components of $I(f)$ are unbounded, and they consist of
\begin{itemize}
 \item[(i)] $I_1$, which also contains any escaping {\wand}s and
any Baker domains of $f$ with at least one boundary point in
$I(f)$,
 \item[(ii)] any Baker domains of $f$ with no boundary points in
$I(f)$ and the infinitely many preimage components of such
Baker domains.
\end{itemize}
\end{itemize}
\end{theorem}
{\it Remark} \quad If every Baker domain of $f$ has a point of
$I(f)$ on its boundary, then the conclusion of
Theorem~\ref{Thm3.1} can be strengthened to `$I(f)$ is
connected'.
\begin{proof}[Proof of Theorem~\ref{Thm3.1}]
Suppose that $E$ is contained in the union of finitely many
components of $I(f)$, say $I_1,I_2,\ldots,I_m$. Take any $z\in
I(f)\cap J(f)$. Since $J(f)\subset\overline{E}$, there exist
$z_n\in E$ such that $z_n\to z$ as $n\to\infty$. Without loss of
generality all terms of this sequence $(z_n)$ lie in a single
component, $I_j$ say. Since $z\in I(f)$, we have $z\in I_j$.
Hence
\begin{equation}\label{eqn5.1}
I(f)\cap J(f)\subset I_1\cup I_2\cup \cdots\cup I_m.
\end{equation}
We now assume that $I_1,I_2,\ldots,I_m$ is the minimal set of
components of $I(f)$ such that (\ref{eqn5.1}) holds. Then $I_j\cap
J(f)\ne\emptyset$, for $j=1,2,\ldots,m$. Note that if the exceptional set $E(f)$ is non-empty, then
\begin{equation}\label{eqn5.1a}
(I_j\setminus E(f))\cap J(f)\ne \emptyset,\quad\text{for } j=1,2,\ldots, m.
\end{equation}
Indeed, if $E(f)=\{\alpha\}\subset I_j\cap J(f)$, then $\alpha$ is a limit point of the backwards orbit of
any non-exceptional point in $I(f)\cap J(f)$ (this follows from
the blowing up property of $J(f)$) and hence $\alpha$ is the limit of a sequence in $I_i\cap J(f)$, say, by (\ref{eqn5.1}). Thus $i=j$ and so (\ref{eqn5.1a}) holds.

If $m=1$, then $I(f)\cap
J(f)$ is contained in one component of $I(f)$, as required.
If $m>1$, then we can take $z_1\in I_1\cap J(f)$ and an open
disc $D$ centred at~$z_1$ so small that
\begin{equation}\label{eq5.2}
D\cap(I_2\cup\cdots \cup I_m)=\emptyset.
\end{equation}
Consider $I_j$, $j\ge 2$. Then
there exists $N\in\N$ such that $f^N(D)$ meets both $I_1\cap
J(f)$ and $I_j\cap J(f)$, by (\ref{eqn5.1a}) and the blowing up property. Hence there
exist $w_1,w_j\in D$ such that
\[
f^N(w_1)\in I_1\cap J(f)\quad\text{and}\quad f^N(w_j)\in I_j\cap J(f),
\]
so $w_1,w_j\in I_1$ by the backwards invariance of $I(f)\cap
J(f)$ and (\ref{eq5.2}). Thus $f^N(I_1)$ is a connected subset of
$I(f)$ that meets both $I_1$ and $I_j$, which is a contradiction.
Hence $m=1$, so part~(a) holds.

Clearly the component $I_1$ is unbounded. If $U$ is a wandering
domain in $I(f)$, then $\partial U\cap I(f)\ne\emptyset$ by
Theorem~\ref{Thm1.1}, so $\partial U\cap I_1\ne\emptyset$ and
hence $U\subset I_1$. Similarly, if $U$ is a Baker domain whose
boundary meets $I(f)$, then $U\subset I_1$; in particular, this
is true if $U$ is completely invariant, since in this case
$\partial U=J(f)$.

To complete the proof of part~(b) we note that any other
component of $I(f)$ must be a Baker domain that has no points of
$I(f)$ in its boundary, and so is not completely invariant, or it must be
one of its infinitely many preimage Fatou components, all of
which are unbounded.
\end{proof}

Theorem~\ref{Thm3.1} includes the following special cases.

\begin{corollary}\label{Cor3.1}
Let $f$ be a {\tef}.
\begin{itemize}
 \item[(a)] Either $I(f)$ is connected or it has infinitely many components.
 \item[(b)]  If $A(f)$ is contained in the
 union of finitely many components of $I(f)$, then all the components of $I(f)$
are unbounded.
\end{itemize}
\end{corollary}
In part~(b) the set $A(f)$ can be replaced by any
subset of $I(f)$ whose closure contains $J(f)$. Such subsets of $I(f)$ include the sets $L(f)$,
$M(f)$ and $Z(f)$, all of which are defined in~\cite{RS09a}.

The proof of Theorem~\ref{Thm3.1}, and hence of
Corollary~\ref{Cor3.1}, used Theorem~\ref{Thm1.1} together with
the blowing up property of $J(f)$. In fact, we can use just the
blowing up property to prove part~(a) of Corollary~\ref{Cor3.1}.
More generally, we have the following result which is not directly related to
escaping points.

\begin{theorem}\label{Thm3.2}
Let $f$ be a {\tef}, and suppose that the set~$E$ is completely
invariant under $f$ and that $J(f)=\overline{E\cap J(f)}$. Then
exactly one of the following holds:
\begin{itemize}
 \item[(1)] $E$ is connected;
 \item[(2)] $E$ has exactly two components, one of which is a singleton $\{\alpha\}$, where $\alpha$ is
 a fixed point of $f$ and $\alpha\in E(f)\cap F(f)$;
 \item[(3)] $E$ has infinitely many components.
\end{itemize}
\end{theorem}
\begin{proof}
We can assume that $E$ is infinite, since any completely
invariant finite set must be a singleton.

Suppose that cases~(1) and~(3) do not hold. Then $E$ has finitely
many components $E_1,E_2,\ldots,E_m$, say, with $m>1$. Since
$J(f)=\overline{E\cap J(f)}$, we can assume that $E_1\ne E(f)$
and that there exists $z_1\in (E_1\setminus E(f))\cap J(f)$. Now
suppose that, for some $j\ge 2$, we have $E_j\ne E(f)$. Choose an
open disc $D$ centred at~$z_1$ and not meeting $E_2\cup\cdots\cup
E_m$. For some $N\in\N$, $f^N(D)$ meets both $E_1$ and $E_j$, by
the blowing up property, so $f^N(E_1)$ meets both $E_1$ and $E_j$
by the backwards invariance of $E$. But $f^N(E_1)$ is a connected
subset of $E$, by the forwards invariance of~$E$, which is a
contradiction.

Thus $E$ has two components $E_1$ and $E(f)=\{\alpha\}$, say.
Then $\alpha\in F(f)$, since any point in $E\cap J(f)$ lies in $E_1$ because it is a limit
point of the backwards orbit of~$z_1$ (this backwards orbit lies in $E$ by the
backwards invariance of $E$).

Finally we show that $f(\alpha)=\alpha$. If
$f(\alpha)\ne\alpha$, then $f(\alpha)\in E_1$, by the complete
invariance of $E$. Since $E_1$ is infinite, $f(\alpha)$ is a limit
point of $E$ and so therefore is $\alpha$, by the complete
invariance of $E$ and the local behaviour of $f$ near~$\alpha$.
Thus we have a contradiction. Hence $\alpha$
is a fixed point of $f$ and case~(2) holds.
\end{proof}
{\it Remarks}

1. The function $f(z)=\tfrac12 z^2e^{2-z}$ shows that case~(2) of
Theorem~\ref{Thm3.2} can occur, with $E=J(f)\cup\{0\}$. In this
case $J(f)$ is connected and $0$ is a fixed point of~$f$ in
$E(f)\cap F(f)$ (see \cite[proof of Theorem~4]{mK98}).

2. Taking $E$ to be the grand orbit of a non-exceptional point in
$J(f)$ shows that in case~(3) of Theorem~\ref{Thm3.2} the set $E$
may have only countably many components.

Note that case~(2) of Theorem~\ref{Thm3.2} cannot occur if we also know,
for example, that $E\subset I(f)$ or that $E$ contains a
neighbourhood of every fixed point of $f$ in $F(f)$. Thus if $E$
is $I(f), L(f), M(f), Z(f)$ or $A(f)$, or the complement of one
of these sets, then $E$ is either connected or has infinitely
many components.

%
%

It is known that if $f$ is a {\tef}, then $J(f)$ is either
connected or it has uncountably many components (see
\cite[Theorem~B]{BD00}), and it is natural to ask if $I(f)$ and
$A(f)$ also have this property.

In Theorem~\ref{Thm3.2} the case when $E$ has infinitely many
components can be strengthened to state that each neighbourhood
of each point $z_0\in J(f)$ must meet infinitely many components
of $E$. This is implied by the following simple result.

\begin{theorem}\label{Thm3.3}
Let $A$ and $B$ be completely invariant sets for a {\tef}~$f$. If, for some $z_0\in J(f)$, there is a neighbourhood
$D$ of $z_0$ such that $A\cap D$ meets only finitely
many components of $B$, then $A$ meets only finitely many
components of $B$.
\end{theorem}
\begin{proof}
Suppose that $A\cap D\subset B_1\cup\cdots\cup B_m$, a
minimal union of components of~$B$. Let $K$ be any compact set,
not meeting the exceptional set $E(f)$. Then, by the blowing up property of $J(f)$, there
exists $N\in\N$ such that $f^N(D)\supset K$. Thus $A\cap K$
meets at most~$m$ components of $B$, by the complete invariance
of $A\cap B$. Hence $A$ meets exactly~$m$ components of $B$ in
$\C\setminus E(f)$, namely $B_1,\ldots, B_m$, and so $A$ meets at
most $m+1$ components of $B$ in $\C$.
\end{proof}
By applying Theorem~\ref{Thm3.3} with $A=\C$ and $B=E$, we deduce
that if $E$ is completely invariant and has infinitely many
components, then each neighbourhood of each point $z_0\in J(f)$
must meet infinitely many components of $E$; in particular, this
is the case in Theorem~\ref{Thm3.2} case~(3).

Other choices of the sets $A$ and $B$ in Theorem~\ref{Thm3.3}
give other corollaries. For example, with $A=A(f)$ and $B=I(f)$,
Theorem~\ref{Thm3.3} and Corollary~\ref{Cor3.1} part~(b),
combine to give the following result.

\begin{corollary}\label{Cor3.3}
Let $f$ be a {\tef}. If, for some $z_0\in J(f)$, there is a
neighbourhood $D$ of $z_0$ such that $A(f)\cap D$ is contained in
the union of finitely many components of $I(f)$, then $A(f)$ is
contained in the union of finitely many components of $I(f)$, so
all the components of $I(f)$ are unbounded.
\end{corollary}

\end{document}